\documentclass[12pt]{article}
\usepackage{amssymb,latexsym,amsmath,euscript}\usepackage{amsfonts}
\usepackage[cp1251]{inputenc}\usepackage[epsf]{}\usepackage[english,russian]{babel}\usepackage{graphicx}
\newtheorem{theorem}{Theorem}
\newtheorem{lemma}{Lemma}[section]
\newtheorem{proposition}{Proposition}
\newtheorem{cor}{Corollary}[section]
\begin{document}

\title{On a classification of axiom A diffeomorphisms with codimension one basic sets and isolated saddles
\footnote{2000{\it Mathematics Subject Classification}. Primary 37D15; Secondary 58C30}
\footnote{{\it Key words and phrases}: axiom A diffeomorphism, basic set, expanding attractor, conjugacy}}

\author{V.~Medvedev$^{1}$\and E.~Zhuzhoma$^{1}$}
\date{}
\maketitle

{\small $^{1}$ National Research University Higher School of Economics, 25/12 Bolshaya Pecherskaya, 603005, \\ Nizhny Novgorod, Russia}

\renewcommand{\abstractname}{Absrtact}
\renewcommand{\refname}{Bibliography}
\renewcommand{\figurename}{Fig.}
%%%%%%%%%%%%%%%%%%%%%%%%%%%%%%%%%%%%

\begin{abstract}
Let $M^n$, $n\geq 3$, be a closed orientable $n$-manifold and $\mathbb{D}_k(M^n;a,b,c)$ the set of axiom A diffeomorp\-hisms $f: M^n\to M^n$ satisfying the following conditions:
(1) $f$ has $k\geq 1$ nontrivial basic sets each is either an orientable codimension one expanding attractor or an orientable codimension one contracting repeller, and other trivial basic sets which are $a$ sinks, $b$ sources, $c$ saddles; (2) the invariant manifolds of isolated saddles are intersected transversally. We classify the diffeomorphisms from $\mathbb{D}_k(M^n;a,b,c)$ up to the global conjugacy on non-wandering sets for the following subsets $\mathbb{S}_k(M^n;a,b,c), \mathbb{P}_k(M^n;0,0,1), \mathbb{M}_k(M^n;0,0,1)$ of $\mathbb{D}_k(M^n;a,b,c)$ where $\mathbb{S}_k(M^n;a,b,c)$ satisfies to the following conditions:

($1_{\mathbb{S}}$) every nontrivial basic set of any $f\in\mathbb{S}_k(M^n;a,b,c)$ is uniquely bunched, and there is at least one nontrivial attractor and at least one nontrivial repeller, i.e. $k\geq 2$;

($2_{\mathbb{S}}$) $c\geq 1$ and all isolated saddles have the same Morse index belonging to $\{1,n-1\}$.

The subset $\mathbb{P}_k(M^n;0,0,1)\subset\mathbb{D}_k(M^n;0,0,1)$ satisfies to the following conditions:

($1_{\mathbb{P}}$) any boundary point of $f\in\mathbb{P}_k(M^n;0,0,1)$ is fixed;

($2_{\mathbb{P}}$) a unique isolated saddle has Morse index different from $\{1,n-1\}$.

The subset $\mathbb{M}_k(M^n;0,0,1)\subset\mathbb{D}_k(M^n;0,0,1)$ satisfies to the following conditions:

($1_{\mathbb{M}}$) any boundary point of $f\in\mathbb{M}_k(M^n;0,0,1)$ is fixed;

($2_{\mathbb{M}}$) a unique isolated saddle has Morse index belonging to $\{1,n-1\}$.

The classification is based on a description of topological structure of supporting manifolds $M^n$.
\end{abstract}

%%%%%%%%%%%%%%%%%%%%%%%%%%%%%%%%%%%%%%%%%%%%%%%%%%%%%%%%%%%%%%%%%%%%%%%%%%%%
\section*{Introduction}\label{s:title-introd}\nopagebreak
%%%%%%%%%%%%%%%%%%%%%%%%%%%%%%%%%%%%%%%%%%%%%%%%%%%%%%%%%%%%%%%%%%%%%%%%%%%%

Dynamical systems satisfying the axiom A (in short, A-systems) was introduced by Smale \cite{Smale67}. By definition, a non-wandering set of A-system has a hyperbolic structure and periodic orbits are dense in the non-wandering set (see basic notation of the Theory of Dynamical Systems in the books \cite{grinmedpoch-book-2016-eng,GrinesZh-book-2021,Robinson-book-99}). According to Mane \cite{Mane88b} and Robinson \cite{Rob76}, A-systems form a wide class containing all structurally stable dynamical systems including Morse-Smale systems and Anosov systems. Here, we consider discrete time A-systems i.e., A-diffeomorphisms $f: M^n\to M^n$ of smooth closed orientable $n$-manifolds, $n\geq 3$.

Due to Smale's Spectral Decomposition Theorem \cite{Robinson-book-99,Smale67}, a non-wandering set of A-diffeomor\-p\-hism splits into pairwise disjoint transitive, closed and invariant pieces called \textit{basic sets}.
A basic set is called \textit{trivial} if it is an isolated periodic orbit. Otherwise, a basic set is \textit{nontrivial}. Bowen \cite{Bow71} proved that the restriction of A-diffeomorphism on a nontrivial basic set has a positive topological entropy i.e. dynamics on nontrivial basic sets are chaotic.

Good examples of nontrivial basic sets are expanding attractors and contracting repellers introduced by Williams \cite{Williams1970a,Williams1974}. He proved that an expanding attractor is locally homeomorphic to the product of Cantor set and Euclidean space. Roughly speaking, the Cantor set is locally the intersection of expanding attractor with stable manifolds while the expanding attractor is the union of unstable manifolds of its points (this was the reason to call such attractors expanding). Many exciting examples of expanding attractors was constructed by Farrel and Jones \cite{FarrelJones1980,FarrelJones1981,Jones1983,Jones1986}. Note that the first examples of expanding attractors, so-called Smale solenoid and DA-attractor, with no names were constructed by Smale \cite{Smale67}. Recall that a DA-attractor, say $\Lambda_a$, is a unique codimension one expanding attractor of A-diffeomorphism $f: \mathbb{T}^n\to\mathbb{T}^n$ of $n$-torus $\mathbb{T}^n$, $n\geq 2$. Such $f$ can be obtained by Smale's surgery from a codimension one Anosov diffeomorphisms $A: \mathbb{T}^n\to\mathbb{T}^n$ such that $f_*=A_*: H_1(\mathbb{T}^n)\to H_1(\mathbb{T}^n)$, see details in \cite{Robinson-book-99}. According to Smale's construction, $\Lambda_a$ is an orientable expanding attractor, and $\Lambda_a=\mathbb{T}^n\setminus W^u(p)$ where $W^u(p)$ is the unstable manifold of a unique isolated source $p$.

Codimension one expanding attractors were completely classified mainly by V. Grines, R. Plykin, Yu. Zhirov, and E. Zhuzhoma \cite{Gr75-77,GinesZhuzhoma1979,GinesZhuzhoma2005,Plykin84}, see the surveys  \cite{GrinesPochinkaZh2014,GrinesZh2006} and the books \cite{grinmedpoch-book-2016-eng,GrinesZh-book-2021}.

The most famous A-diffeomorphisms with all trivial basic sets are Morse-Smale diffeo\-mor\-p\-hisms \cite{Palis69}. Recall that Morse-Smale diffeomorphisms are structurally stable ones with zero topological entropy. Therefore they are satisfied the so-called strong transversality condition, i.e. invariant manifolds of saddle periodic orbits are intersected transversally.
Last decades, the progress in the problem of classification of Morse-Smale diffeomorphisms was obtained due to
Ch. Bonatti, V. Grines, E. Gurevich, V. Medvedev, and O. Pochinka, see the surveys \cite{GrinesGurevichZhPochinka2019,MedvedevZhuzhoma2008-MIAN} and the books \cite{GrinGurevich-book-2022,grinmedpoch-book-2016-eng}.

It is natural to consider the problem of classification for A-diffeomorphisms whose nontrivial basic sets are codimension one expanding attractors and contracting repellers while trivial basic sets demonstrate Morse-Smale dynamics. Following \cite{GrinMedvZh-2022,GrinMedvZh-2024} we introduce the set $\mathbb{D}_k(M^n;a,b,c)$ of A-diffeomorp\-hisms $f: M^n\to M^n$ satisfying the following conditions:

($1_{\mathbb{D}}$) every $f\in\mathbb{D}_k(M^n;a,b,c)$ has $k\geq 1$ nontrivial basic sets each is either an orientable codimension one expanding attractor or an orientable codimension one contracting repeller, and other basic sets are trivial ones which are $a$ sinks, $b$ sources, and $c$ saddles;

($2_{\mathbb{D}}$) the invariant manifolds of isolated saddles are intersected transversally.

The set $\mathbb{D}_k(M^n;a,b,c)$ is obviously non-empty. For example, the DA-diffeomorphism mentioned above belongs to $\mathbb{D}_k(\mathbb{T}^n;0,1,0)$. Robinson and Williams \cite{RobinsonWilliams71} constructed a diffeomor\-p\-hism $f\in\mathbb{D}_2(\mathbb{T}^n\sharp\mathbb{T}^n;0,0,0)$ where $\mathbb{T}^n\sharp\mathbb{T}^n$ is a connected sum of two torii $\mathbb{T}^n$. For Reader's convenience, let us explain this construction. Take a DA-diffeomorphism $f_0: \mathbb{T}^n\to\mathbb{T}^n$ those non-wandering set consists of an isolated source $s_0$ and codimension one orientable expanding attractor $\Lambda$. Then the diffeomorphism $f_0^{-1}$ has the non-wandering set consisting of the sink $s_0$ denoted by $s_1$ and the codimension one orientable contracting repeller $\Lambda$ denoted by $\Lambda_1$. We can assume that $f_0^{-1}$ is defined on a copy of $\mathbb{T}^n$. Deleting small neighborhoods of $s_0$ and $s_1$ one can construct a connected sum $M^n=\mathbb{T}^n\sharp\mathbb{T}^n$ on which $f_0$ and $f_0^{-1}$ induce $f\in\mathbb{D}_2(\mathbb{T}^n\sharp\mathbb{T}^n;0,0,0)$ those non-wandering set consists of the orientable codimension one expanding attractor $\Lambda$ and contracting repeller $\Lambda_1$.

Our goal is a classification of some classes of diffeomorphisms from the set $\mathbb{D}_k(M^n;a,b,c)$.
There are several approaches to classify dynamical systems and invariant sets, and a classification depends on the type of conjugacy.
Recall that two maps $f: M\to M$, $g: N\to N$ are (topologically) \textit{conjugate} provided there is a homeomorphism $h: M\to N$ such that $h\circ f=g\circ h$.
It is a difficult problem to classify dynamical systems under conjugacy mappings on whole supporting manifold. The first natural step is a classification of restrictions of diffeomorphisms on special invariant subsets.
The second natural step is to ask, when two diffeomorphisms are conjugate in neighborhoods of their invariant sets? Robinson and Williams \cite{RobinsonWilliams76} constructed two diffeomorphisms $f$ and $g$ of non-homeomorphic 5-manifolds with expanding 2-dimensional attractors $\Lambda _f$ and $\Lambda _g$ respectively such that the restriction $f|_{\Lambda_f}: \Lambda_f\to\Lambda_f$ is conjugate to the restriction
$g|_{\Lambda_g}: \Lambda_g\to\Lambda_g$ but there is not even a homeomorphism from a neighborhood of $\Lambda_f$ to a neighborhood of $\Lambda_g$ taking $\Lambda_f$ to $\Lambda_g$. Another examples see in \cite{IsaenkovaZh2009}, where the first type of conjugacy (i.e. a conjugacy of restrictions) is called an intrinsic conjugacy while the second type of conjugacy (when a conjugacy map is defined on a neighborhood of invariant set) is called a neighbor conjugacy. Clearly, a neighbor conjugacy implies an intrinsic conjugacy, and the neighbor conjugacy takes into account embedding of invariant sets in supporting manifolds. Here, we consider a global conjugacy which can be considered as an intermediate type. To be precise, suppose diffeomorphisms $f,f': M^n\to M^n$ have invariant sets $\Omega$, $\Omega^{\prime}$ respectively. We say that $f$, $f^{\prime}$ are
\textit{globally conjugate on the sets} $\Omega$, $\Omega^{\prime}$ if there is a homeomorphism $h: M^n\to M^n$ such that
$h(\Omega) = \Omega^{\prime}$ and $f^{\prime}\vert_{\Omega^{\prime}} = h\circ f\circ h^{-1}\vert_{\Omega^{\prime}}.$

To formulate the main results we introduce some notation. Suppose $\Lambda$ is an orientable codimension one expanding attractor (similar notation holds for a contracting repeller) of an A-diffeomorphism $f$. Then any stable manifold $W^s(x)$, $x\in\Lambda$, is one-dimensional and $W^s(x)\backslash x$ consists of two components. Due to \cite{Gr75-77} for $n=2$, and \cite{GinesZhuzhoma2005}, Lemmas 1.2, 1.5 for $n\geq 3$, at least one component of $W^s(x)\backslash x$ intersects $\Lambda$. A point $x\in\Lambda$ is called \textit{boundary} if there is a component of $W^s(x)\backslash x$ denoted by $W_{\emptyset}^{s}(x)$ which does not intersect $\Lambda$. The set $B(f)\subset\Lambda$ of boundary points is non-empty, invariant and finite (so every boundary point is periodic). Let $p_1, \ldots, p_r\in B(f)$ be all boundary points such that $W^{s}_{\emptyset}(p_1)$, $\ldots$, $W^{s}_{\emptyset}(p_r)$ belong to the same component of $W^s({\Lambda})\setminus\Lambda$. The union $\cup_{i=1}^rW^u(p_i)$ denoted by $b^u$ is called a \textit{bunch}, $r$ is a \textit{degree} of the bunch $b^u$, and $p_1, \ldots, p_r$ are called \textit{associated} points. Since $\Lambda$ is orientable, each bunch of $\Lambda$ has degree two \cite{Gr75-77,Plykin74,Plykin84}, and hence a bunch has two associated periodic points \cite{GinesZhuzhoma2005} (see also the books \cite{grinmedpoch-book-2016-eng,GrinesZh-book-2021}). We call $\Lambda$ \textit{uniquely bunched} if $\Lambda$ has a unique bunch. In this case the corresponding associated periodic points are fixed.

We consider three following subsets
 $$\mathbb{S}_k(M^n;a,b,c),\quad \mathbb{P}_k(M^n;a,b,c),\quad \mathbb{M}_k(M^n;a,b,c)\subset\mathbb{D}_k(M^n;a,b,c)$$
where the first subset $\mathbb{S}_k(M^n;a,b,c)\subset\mathbb{D}_k(M^n;a,b,c)$ satisfies to the following conditions:

($1_{\mathbb{S}}$) every nontrivial basic set of any $f\in\mathbb{S}_k(M^n;a,b,c)$ is uniquely bunched, and there is at least one nontrivial attractor and at least one nontrivial repeller, i.e. $k\geq 2$;

($2_{\mathbb{S}}$) $c\geq 1$ and all isolated saddles have the same Morse index belonging to $\{1,n-1\}$.

\medskip
The second subset $\mathbb{P}_k(M^n;0,0,1)\subset\mathbb{D}_k(M^n;0,0,1)$ satisfies to the following conditions:

($1_{\mathbb{P}}$) any boundary point of $f\in\mathbb{P}_k(M^n;0,0,1)$ is fixed;

($2_{\mathbb{P}}$) a unique isolated saddle has Morse index different from $\{1,n-1\}$.

\medskip
The third subset $\mathbb{M}_k(M^n;0,0,1)\subset\mathbb{D}_k(M^n;0,0,1)$ satisfies to the following conditions:

($1_{\mathbb{M}}$) any boundary point of $f\in\mathbb{M}_k(M^n;0,0,1)$ is fixed;

($2_{\mathbb{M}}$) a unique isolated saddle has Morse index belonging to $\{1,n-1\}$.

\medskip
A crucial role in our classification plays a complete invariant of conjugacy for a codimension one basic set constructed by R. Plykin and V. Grines. Based on his invariant, we construct invariants of conjugacy for a diffeomorphism $f\in\mathbb{S}_k(M^n;a,b,c)\cup\mathbb{P}_k(M^n;0,0,1)\cup\mathbb{M}_k(M^n;0,0,1)$.

For $f\in\mathbb{S}_k(M^n;a,b,c)$, one defines a $k$-tube $t(f)$ which is a family $(f_1,\ldots,f_k)$ of codimension one Anosov automorphisms $f_i: \mathbb{T}^n\to\mathbb{T}^n$, see details in Section \ref{s:prev}. We introduce a set $A_k$ of admissible $k$-tubes and the notion of equivalence for $k$-tubes as well. Due to condition ($2_{\mathbb{S}}$), one can define $ind (f)$ to be Morse index of an (any) isolated saddle of $f$.

Our constructions are intimately connected with a topological structure of supporting manifolds. Therefore we include the description of topological structure of supporting manifolds in the main results.

Below, $\mathbb{S}^m$ is an $m$-sphere, $\mathbb{T}^n$ is an $n$-torus.
\begin{theorem}\label{thm:classif-for-S}
Diffeomorphisms $f_i\in\mathbb{S}_k(M^n;a_i,b_i,c_i)$, $i=1,2$, are globally conjugate on their non-wandering sets if and only if the following conditions hold:
\begin{itemize}
  \item $k$-tubes $t(f_1)$, $t(f_2)$ are equivalent under linear automorphisms $\zeta_1$, $\ldots$, $\zeta_k: \mathbb{T}^n\to\mathbb{T}^n$ such that the determinants of all $\zeta_j$ have the same sign;
  \item $a_1+b_1=a_2+b_2$, $c_1=c_2$.
\end{itemize}
Besides,
\begin{equation}\label{eq:S-decomposition-manifold}
    M^n=\underbrace{\mathbb{T}^n\sharp\cdots\sharp\mathbb{T}^n}_{k\geq 2}.
\end{equation}
Moreover, if $f\in\mathbb{S}_k(M^n;a,b,c)$ then $k+a+b=c+2$ and the $k$-tube $k(f)$ is admissible, and $t(f)$ agreed with the triple $(a,b,c)$.
Conversely, given any $k$-tube $t\in A_k$ and integers $k\geq 2$, $a,b\geq 0$, $c\geq 1$ with $k+a+b=c+2$ such that $t$ agreed with the triple $(a,b,c)$, there is a manifold $M^n$ of the kind (\ref{eq:S-decomposition-manifold}) and a diffeomorphism $f\in\mathbb{S}_k(M^n;a,b,c)$ with $k(f)=t$.
\end{theorem}

For the set $\mathbb{P}_k(M^n;0,0,1)$, it is easy to see that if $\mathbb{P}_k(M^n;0,0,1)\neq\emptyset$ then $k\geq 2$ (there are at least one non-trivial attractor and at least one non-trivial repeller). Moreover, one can prove that the dimension $n\in\{4,8,16\}$. Given any $k\geq 2$ and $n\in\{8,16\}$, we construct an invariant of global conjugacy for every $f\in\mathbb{P}_k(M^n;0,0,1)$ which is a graph $\Gamma_{\mathbb{P}}(f)$ endowing with an additional information, and one introduces the definition of commensurability of graphs. In addition, we introduce the set $\Gamma^k_{\mathbb{P}}$, $k\geq 2$, of abstract graphs (see details in Section \ref{s:prev}).
The following result says that the graph $\Gamma_{\mathbb{P}}(f)$ of $f\in\mathbb{P}_k(M^n;0,0,1)$, $k\geq 2$, $n\in\{8,16\}$, up to a commensurability is a complete invariant of global conjugacy. It is convenient to assume $\underbrace{K^n\sharp\cdots\sharp K^n}_{g\geq 0}=\mathbb{S}^n$ provided $g=0$.
\begin{theorem}\label{thm:classif-for-P}
Two diffeomorphisms $f_1,f_2\in\mathbb{P}_k(M^n;0,0,1)$, $k\geq 2$ and $n\in\{8,16\}$, are globally conjugate on their non-wandering sets if and only if the graphs $\Gamma_{\mathbb{P}}(f_1)$, $\Gamma_{\mathbb{P}}(f_2)$ are commensurable. Any graph $\Gamma_{\mathbb{P}}(f)$ belongs to $\Gamma^k_{\mathbb{P}}$.
Given any graph $\gamma\in\Gamma^k_{\mathbb{P}}$, there are a closed smooth connected orientable $n$-manifold $M^n$, $n\in\{8,16\}$, and a diffeomorphism $f\in\mathbb{P}_k(M^n;0,0,1)$ such that $\gamma=\Gamma_{\mathbb{P}}(f)$. In addition, a supporting manifold $M^n$ for any $f\in\mathbb{P}_k(M^n;0,0,1)$, $k\geq 2$, is homeomorphic to
\begin{equation}\label{eq:P-decomposition-structure}
    M^n=\underbrace{\mathbb{T}^n\sharp\cdots\sharp\mathbb{T}^n}_{k\geq 2}\sharp\underbrace{\left({S}^{n-1}\times {S}^1\right)\sharp\cdots\sharp\left({S}^{n-1}\times {S}^1\right)}_{g\geq 0}\sharp N^n
\end{equation}
where
\begin{itemize}
  \item $n\in\{4,8,16\}$
  \item $N^n$ is a simply connected manifold homeomorphic to $\mathbb{S}^{\frac{n}{2}}\sqcup\mathbb{B}^n$, here $\mathbb{B}^n$ is an open $n$-ball;
  \item $N^n$ is a projective-like manifold provided $n\in\{8,16\}$.
\end{itemize}
\end{theorem}

Similar result holds for $\mathbb{M}_k(M^n;0,0,1)$. Here for $f\in\mathbb{M}_k(M^n;0,0,1)$, we construct an invariant of global conjugacy which is a graph $\Gamma_{\mathbb{M}}(f)$ endowing with an additional information. One introduces the notion of commensurability, and the set $\Gamma^k_{\mathbb{M}}$, $k\geq 2$, of abstract graphs.
\begin{theorem}\label{thm:classif-for-M}
Two diffeomorphisms $f_1,f_2\in\mathbb{M}_k(M^n;0,0,1)$, $k\geq 2$, are globally conjugate on their non-wandering sets if and only if the graphs $\Gamma_{\mathbb{M}}(f_1)$, $\Gamma_{\mathbb{M}}(f_2)$ are commensurable. Any graph $\Gamma_{\mathbb{M}}(f)$ belongs to $\Gamma^k_{\mathbb{M}}$.
Given any graph $\gamma\in\Gamma^k_{\mathbb{M}}$, there are a closed smooth connected orientable $n$-manifold $M^n$ and a diffeomorphism $f\in\mathbb{M}_k(M^n;0,0,1)$ such that $\gamma=\Gamma_{\mathbb{M}}(f)$. In addition, a supporting manifold $M^n$ for any $f\in\mathbb{P}_k(M^n;0,0,1)$, $k\geq 2$, is homeomorphic to
 $$ M^n=\underbrace{\mathbb{T}^n\sharp\cdots\sharp\mathbb{T}^n}_{k\geq 2}\sharp\underbrace{\left({S}^{n-1}\times {S}^1\right)\sharp\cdots\sharp\left({S}^{n-1}\times {S}^1\right)}_{g\geq 0}. $$
\end{theorem}

Let us mention some results concerning the subject of the paper. The classification with no trivial basic sets was considered in \cite{GrinMedvZh-2024}.
In \cite{MedvedevZh2023}, it was introduced the following four types of A-diffeomorphisms: regular, semi-chaotic, chaotic, and super chaotic ones (to be precise, such types were introduced for a wide class of Smale A-homeomorphisms). A family of basic sets of A-diffeomorphisms $f$ is naturally divided into sink basic sets $\omega(f)$, source basic sets $\alpha(f)$, and saddle basic sets $\sigma(f)$. We say that $f$ is \textit{regular} if all basic sets $\omega(f)$, $\sigma(f)$, $\alpha(f)$ are trivial while $f$ is \textit{semi-chaotic} if exactly one family from the families $\omega(f)$, $\sigma(f)$, $\alpha(f)$ consists of non-trivial basic sets, and $f$ is  \textit{chaotic} if exactly two families from the families $\omega(f)$, $\sigma(f)$, $\alpha(f)$ consists of non-trivial basic sets, and at last $f$ is \textit{super chaotic} if all basic sets $\omega(f)$, $\sigma(f)$, $\alpha(f)$ are non-trivial. In \cite{MedvedevZh2023}, it was formulated necessary and sufficient condi\-ti\-ons of conjugacy for regular, semi-chaotic, and chaotic A-diffeomorphisms provided chaotic A-diffeomorphisms have either trivial sink basic sets or trivial source basic sets. We see that the A-diffeomorphisms under consideration belongs to the set of chaotic A-diffeomorphisms, but each of them has non-trivial sink and source basic sets. Thus, the main result of \cite{MedvedevZh2023} does not cover the main results of our paper.

The structure of the paper is the following. In Section \ref{s:prev}, we formulate basic definitions and previous results. In section \ref{s:topol}, one describes topological structures of supporting manifolds. In Section \ref{s:proofs}, we prove main results.

\textsl{Acknowledgments}. This work is an output of a research project implemented as part of the Basic Research Program at the National Research University Higher School of Economics.

%%%%%%%%%%%%%%%%%%%%%%%%%%%%%%%%%%%%%%%%%%%%%%%%%%%%%%%%%%%%
\section{Basic definitions and previous results}\label{s:prev}
%%%%%%%%%%%%%%%%%%%%%%%%%%%%%%%%%%%%%%%%%%%%%%%%%%%%%%

\textsl{A-diffeomorphisms}.
$f: M^n\to M^n$ is called an \textit{A-diffeomorphism} if its non-wandering set $NW(f)$ is hyperbolic and periodic points are dense in $NW(f)$ \cite{Smale67}. The existence of hyperbolic structure implies the existence of unstable $W^u(x)$ and stable $W^s(x)$ manifolds respectively for every point $x\in NW(f)$ \cite{HP70}.
Due to Smale's Spectral Decomposition Theorem, the non-wandering set $NW(f)$ is a finite union of pairwise disjoint $f$-invariant closed sets $\Omega_1$, $\ldots$, $\Omega_k$
such that every restriction $f|_{\Omega_i}$ is topologically transitive. These $\Omega_i$ are called \textit{basic sets}. The dimension $\dim W^u(x)$, $x\in\Omega_i$, is called a \textit{Morse index} of $\Omega_i$.
A basic set $\Omega$ is an \textit{expanding attractor} if $\Omega$ is an attractor and the topological dimension $\dim\Omega$ equals Morse's index of $\Omega$ \cite{Williams1974}. A basic set $\Lambda$ of $f$ is called a \textit{contracting repeller} provided $\Lambda$ is an expanding attractor of $f^{-1}$.

By definition, let $W^s_{\epsilon}(x)\subset W^s(x)$ (resp. $W^u_{\epsilon}(x)\subset W^u(x))$ be the $\epsilon$-neighborhood of $x$ in the intrinsic topology of the manifold $W^s(x)$ (resp. $W^u(x))$, where $\epsilon >0$. We say that a basic set $\Omega$ is \textit{orientable} provided the index of intersection $W^s_{\alpha}(x)\cap W^u_{\beta}(x)$ is the same at each point of this intersection for any $\alpha>0$, $\beta>0$, $x\in\Omega$. Smale's solenoid is an orientable attractor \cite{Smale67} while Plykin attractor is non-orientable \cite{Plykin74}.

\medskip
\textsl{Codimension one basic sets}.
Any codimension one expanding attractor $\Lambda$ consists of $(n-1)$-dimensional unstable manifolds $W^u(x)$, $x\in\Lambda$ and locally homeomorphic to the product of $(n-1)$-dimensional Euclidean space and Cantor set of unit interval \cite{Plykin71}, \cite{Williams1974}. An unstable manifold $W^u_x$ passing through a boundary point of the Cantor set is called \textit{boundary}. One can prove that there are only finitely many boundary unstable manifolds each passes through a \textit{boundary} periodic point. The boundary unstable manifolds of $\Lambda$ split into a finite number of bunches \cite{Gr75-77}, \cite{GinesZhuzhoma2005}, \cite{Plykin74}.
If $\Omega$ is a codimension one orientable expanding attractor of A-diffeomorphism $f: M^n\to M^n$, $n\geq 3$, then $\Omega$ has only 2-bunches \cite{Plykin84}. Suppose $b^u_{pq}=W^u(p)\cup W^u(q)$ is a 2-bunch where $p$, $q$ are boundary periodic points of $\Omega$. The points $p$, $q$ are called \textit{associated}, and there is a one-to-one correspondence between bunches and pairs of associated points.

\medskip
\textsl{Canonical neighborhoods}.
Let $\Omega$ be a codimension one connected orientable expanding attractor or contracting repeller of A-diffeomorphism $f$. Plykin \cite{Plykin84} proved that $\Omega$ has a neighborhood $U(\Omega)$ which is homeomorphic to  $\mathbb{T}^n$ with $l$ deleted closed $n$-disks $D^n_i$, i.e. $U(\Omega)=\mathbb{T}^n\setminus\cup_{i=1}^lD^n_i$. Moreover, either $U(\Omega)$ is an attracting neighborhood for $f$ provided $\Omega$ is an attractor or $U(\Omega)$ is an attracting neighborhood for $f^{-1}$ provided $\Omega$ is a repeller. Such $U(\Omega)$ is called a \textit{canonical neighborhood}.

The next result follows from \cite{GinesZhuzhoma2005}, Corollary 2.1 and Lemma 3.1 (see also Theorem 5.1). Recall that a DA-diffeomorphism $g: \mathbb{T}^n\to\mathbb{T}^n$ is an A-diffeomorphism such that either $NW(g)$ consists of codimension one orientable expanding attractor and finitely many isolated source periodic orbits or $NW(g)$ consists of codimension one orientable contracting repeller and finitely many isolated sink periodic orbits. Note that a nontrivial basic set of DA-diffeomorphism is automatically connected \cite{GinesZhuzhoma2005}. A DA-diffeomorphism is \textit{simplest} if its non-trivial basic set is uniquely bunched.
\begin{lemma}\label{lm:nbhd-of-attractor-omega}
Let $\Omega$ be a codimension one connected orientable expanding attractor or contracting repeller with $l$ bunches of A-diffeomor\-p\-hism $f: M^n\to M^n$, $n\geq 3$, and $U(\Omega)$ a canonical neighborhood of $\Omega$. Then the restriction $f|_{U(\Omega)}: U(\Omega)\to f(U(\Omega))$ is extended to a DA-diffeomorp\-hism $\tilde{f}: \mathbb{T}^n\to\mathbb{T}^n$, and $\tilde{f}$ is defined up to a homotopy trivial conjugacy. In addition, there is $r\in\mathbb{N}$ such that the spheres $S_i$, $f^r(S_i)$ bound a domain homeomorphic to $S^{n-1}\times [0;1]$ for all $i=1,\ldots,l$ where $S_i$ are boundary components of $U(\Omega)$. If any boundary point of $\Omega$ is fixed then $r=1$. Moreover, if $\Omega$ is uniquely bunched, then $U(\Omega)$ has a unique boundary component and $\tilde{f}$ is a simplest DA-diffeomorphism.
\end{lemma}

\medskip
\textsl{Plykin-Grines invariant of conjugacy}.
Suppose $\Omega$ and $f: M^n\to M^n$, $n\geq 3$, satisfy the condition of Lemma \ref{lm:nbhd-of-attractor-omega}. Let $\tilde{f}: \mathbb{T}^n\to\mathbb{T}^n$ be a corresponding DA-diffeomorphism. Then $\tilde{f}$ induces the linear isomorphism $\tilde{f}_*: H_1(\mathbb{T}^n,\mathbb{Z})\to H_1(\mathbb{T}^n,\mathbb{Z})$ which can be considered as the automorphism $\tilde{f}_*: \mathbb{T}^n\to\mathbb{T}^n$. Grines and Zhuzhoma \cite{GinesZhuzhoma1979}, Theorem 1, proved that $\tilde{f}_*$ is a codimension one Anosov automorphism. It follows from \cite{Fr70}, Theorem 2.2, that there is a homotopy trivial continuous map $h: \mathbb{T}^n\to\mathbb{T}^n$ that semi-conjugates $f$ to $\tilde{f}_*$, i.e. $h\circ f=\tilde{f}_*\circ h$. As a consequence, the stable manifolds of $\tilde{f}_*$ is one-dimensional provided $\Omega$ is an attractor, and the unstable manifolds of $\tilde{f}_*$ is one-dimensional provided $\Omega$ is a repeller. According \cite{GinesZhuzhoma1979}, Theorem 1 (see also \cite{GinesZhuzhoma2005}), $h(\Omega)=\mathbb{T}^n$ and $h(B_f)=P_f$ is a finite set of periodic orbits of the Anosov automorphism $\tilde{f}_*$ where $B_f$ is the set of boundary points of $f$. Thus, the basic set $\Omega$ of $f$ corresponds a triple $(\tilde{f}_*,P,\epsilon)$ where $\epsilon=a$ provided $\Omega$ is an attractor and $\epsilon=r$ provided $\Omega$ is a repeller.

Two triples $(\tilde{f}_*,P,\epsilon)$, $(\tilde{f'}_*,P',\epsilon')$ are called \textit{equivalent} provided $\epsilon=\epsilon'$, and there is an automor\-p\-hism $\zeta: \mathbb{T}^n\to\mathbb{T}^n$ such that $\tilde{f'}_*=\zeta\circ\tilde{f}_*\circ\zeta^{-1}$ and $\zeta(P)= P'$. Grines \cite{Gr75-77} for $n=2$ and Plykin \cite{Plykin84} for $n\geq 3$ proved that the triple $(\tilde{f}_*,P,\epsilon)$ is a complete invariant of global conjugacy for $f|_{\Omega}$ up to the equivalence. For references, we formulate this classical result.
\begin{theorem}\label{thm:plyk-gr-class}
Let $(A,P,\epsilon)$ be a triple where $A: \mathbb{T}^n\to\mathbb{T}^n$ is a codimension one Anosov automorphism, $P$ is a finite family of periodic orbits of $A$, and either $\epsilon=a$ provided the stable manifolds of $A$ is one-dimensional or $\epsilon=r$ provided the unstable manifolds of $A$ is one-dimensional. Then there is a DA-diffeomorphism $f: \mathbb{T}^n\to\mathbb{T}^n$ with a codimension one orientable basic set $\Omega$ such that

1) $\Omega$ is either an expanding attractor provided $\epsilon=a$ or a contracting repeller provided $\epsilon=r$;

2) $f$ is semi-conjugates to $A$ by a homotopy trivial continuous map $h: \mathbb{T}^n\to\mathbb{T}^n$, that is $h\circ f=A\circ h$;

3) $h(B_f)=P$ where $B_f$ is the set of boundary points of $f$.

Let $f,f: M^n\to M^n$ be A-diffeomorphisms and $\Omega$, $\Omega'$ its codimension one orientable basic sets respectively (each is either an expanding attractor or contracting repeller). Then $f$, $f'$ are globally conjugate on $\Omega$, $\Omega'$ respectively if and only if the triples $(\tilde{f}_*,P,\epsilon)$, $(\tilde{f'}_*,P',\epsilon')$ are equivalent.
\end{theorem}

%******************************************

%As a consequence, one gets the following result for a simplest DA-diffeomorphism.
%\begin{lemma}\label{lm:conjugacy-of-simplest-da}
%Two simplest DA-diffeomorphisms $f_1$, $f_2$ are globally conjugate on its non-trivial basic sets if and only if $f_{1*}$, $f_{2*}$ are conjugate under a linear automorp\-hism $\zeta: \mathbb{T}^n\to\mathbb{T}^n$.
%\end{lemma}

%***************************************************

Let $\varphi: S\to\varphi(S)$ be a homeomorphism and $S$, $\varphi(S)$ the submanifolds of $\mathbb{T}^n$. Then the orientation of $\mathbb{T}^n$ induces interior orientations on $S$ and $\varphi(S)$. One says that
$\varphi: S\to\varphi(S)$ preserves orientation if $\varphi$ preserves the interior orientations of $S$ and $\varphi(S)$ \cite{Hatcher-book}. For references, we formulate the following statement from \cite{GrinMedvZh-2024}, Proposition 1.
\begin{lemma}\label{lm:orient-determinant}
Let $g, g': \mathbb{T}^n\to\mathbb{T}^n$ be DA-diffeomorphisms and $\Omega$, $\Omega'$ nontrivial basic sets of $g$, $g'$ respectively. Suppose a homeomorphism
$\varphi: \mathbb{T}^n\to\mathbb{T}^n$ is a global conjugacy of $g$ and $g'$ on $\Omega$, $\Omega'$ respectively, so that the triples $(A,P,\epsilon)$, $(A',P',\epsilon')$ are equivalent, i.e. there is an automorphism
$\zeta: \mathbb{T}^n\to\mathbb{T}^n$ that conjugates Anosov diffeomorphisms $A$, $A'$, and $\zeta(P)=P'$. Let $S$ be a characteristic sphere of $\Omega$. Then the restriction $\varphi|_{S}: S\to\varphi(S)$ preserves orientation if and only if the determinant of $\zeta$ is positive (hence, $\varphi|_{S}$ reverses orientation if and only if the determinant of $\zeta$ is negative).
\end{lemma}

\medskip
\textsl{Invariant of conjugacy for $\mathbb{S}_k(M^n;a,b,c)$}.
A $k$-\textit{tube of Anosov automorphisms} (or simply, $k$-\textit{tube}) is $t=(A_1,\ldots,A_k)$ where $A_i: \mathbb{T}^n\to\mathbb{T}^n$ is a codimension one Anosov automorphism, $i=1,\ldots,k$. It follows from Theorem \ref{thm:plyk-gr-class} that every $f\in\mathbb{D}_k(M^n;a,b,c)$ corresponds a $k$-tube $t(f)=\left(\tilde{f}_{1*},\ldots,\tilde{f}_{k*}\right)$. Lemma \ref{lm:nbhd-of-attractor-omega} allows to clarify this result for $\mathbb{S}_k(M^n;a,b,c)$ as follows.
\begin{lemma}\label{lm:diffeo-gives-tube}
Suppose $f\in\mathbb{S}_k(M^n;a,b,c)$. Then $f$ corresponds a $k$-tube $t(f)=\left(\tilde{f}_{1*},\ldots,\tilde{f}_{k*}\right)$ consisting of Anosov automorphisms $\tilde{f}_{i*}: \mathbb{T}^n\to\mathbb{T}^n$ induced by simplest DA-diffeomorphisms $f_i$, $i=1,\ldots,k$.
\end{lemma}

Denote by $t_u$ the number of Anosov automorphisms with one-dimensional unstable manifolds.
A $k$-tube $(A_1,\ldots,A_k)$ is \textit{admissible} if $t_u=1$ or $t_u=k-1$. Another words, one of the following possibilities holds :

1) there is a unique Anosov automorphism, say $A_j$, with one-dimensional stable manifolds while another $A_i$, $i\neq j$, have one-dimensional unstable manifolds;

2) there is a unique Anosov automorphism, say $A_j$, with one-dimensional unstable manifolds while another $A_i$, $i\neq j$, have one-dimensional stable manifolds.

We'll say that two $k$-tubes $(A_1,\ldots,A_k)$, $(A_1',\ldots,A_k')$ are \textit{equivalent} provided there are a $k$-permutation $\tau$ and linear automorphisms $\zeta_1$, $\ldots$, $\zeta_k : \mathbb{T}^n\to\mathbb{T}^n$ such that $\zeta_i\circ A_i=A_{\tau(i)}'\circ\zeta_i$ for every $i=1,\ldots,k$. Sometimes we'll say that $(A_1,\ldots,A_k)$, $(A_1',\ldots,A_k')$ are equivalent under linear automorphisms $\zeta_1$, $\ldots$, $\zeta_k$. Clearly that if $k$-tubes $t$, $t'$ are equivalent, then $t_u=t'_u$.

We'll say that the admissible $k$-tube $t$ \textit{agreed} with the triple $(a,b,c)$ where $a,b\geq 0$, $c\geq 1$, $k+a+b=c+2$, when one of the following case holds:
\begin{itemize}
  \item $b=0$ provided $t_u=1$;
  \item $a=0$ provided $t_u=k-1$.
\end{itemize}

We'll show that admissible tubes agreeing with triples $(a,b,c)$ form complete invariants for $\mathbb{S}_k(M^n;a,b,c)$.

\medskip
\textsl{Preliminaries on topological structures}.
A closed manifold $M^n$ is called \textit{projective-like} provided
(1) $n\in\{2,4,8,16\}$; (2) there is a locally flat embedded $\frac{n}{2}$-sphere $S^{\frac{n}{2}}\subset M^n$ such that $M^n\setminus S^{\frac{n}{2}}$ is an open $n$-ball \cite{MedvedevZhuzhoma2016}.

Suppose the boundary $\partial K^n$ of compact $n$-manifold $K^n$ consists of $(n-1)$-spheres. It is convenient to use the notation $\widehat{K^n}$ for manifold obtained from $K^n$ by capping off each $(n-1)$-sphere component of $\partial K^n$ with an $n$-cell. Thus $\widehat{K^n}=K^n\left(\cup_{i=1}^sB^n_i\right)$ is a closed manifold where $s$ is the number of boundary components of $\partial K^n$.

The following result concerns to a topological structure of supporting manifolds and invariants of conjugacy for $f\in\mathbb{P}_k(M^n;0,0,1)\cup\mathbb{M}_k(M^n;0,0,1)$.
\begin{lemma}\label{lm:to-top-str-for-P-and-M}
Suppose $f\in\mathbb{P}_k(M^n;0,0,1)\cup\mathbb{M}_k(M^n;0,0,1)$.
Let $\Omega_1$, $\ldots$, $\Omega_k$ be nontrivial basic sets of $f$, and $U(\Omega_1)$, $\ldots$, $U(\Omega_k)$ its canonical neighborhoods respectively, and $\sigma$ a unique isolated saddle of $f$.
Then the connected components of $M^n\setminus\left(\cup_{i=1}^kU(\Omega_i)\right)$ form a disjoint union $K_{\sigma}\cup_{\nu=1}^rR_{\nu}$ where $K_{\sigma}$ contains $\sigma$ while at every $R_{\nu}$ there are no non-wandering points. Besides, every $R_{\nu}$ is homeomorphic to $\mathbb{S}^{n-1}\times [0;1]$ where one boundary component, say $\mathbb{S}^{n-1}\times\{0\}$, is a boundary component of some canonical neighborhood of expanding attractor while $\mathbb{S}^{n-1}\times\{1\}$ is a boundary component of some canonical neighborhood of contracting repeller.

Moreover, if $f\in\mathbb{P}_k(M^n;0,0,1)$ then $\widehat{K_{\sigma}}=K_{\sigma}\cup B_1^n\cup B_2^n$ and
\begin{itemize}
  \item $n\in\{4,8,16\}$;
  \item $\widehat{K_{\sigma}}$ is a simply connected manifold homeomorphic to $\mathbb{S}^{\frac{n}{2}}\sqcup\mathbb{B}^n$, here $\mathbb{B}^n$ is an open $n$-ball;
  \item $\widehat{K_{\sigma}}$ is a projective-like manifold provided $n\in\{8,16\}$.
\end{itemize}
If $f\in\mathbb{P}_k(M^n;0,0,1)$ then $\widehat{K_{\sigma}}=K_{\sigma}\cup B_1^n\cup B_2^n\cup B_3^n=\mathbb{S}^n$
\end{lemma}
\textsl{Proof}. The equality $M^n\setminus\left(\cup_{i=1}^kU(\Omega_i)\right)=K_{\sigma}\cup_{\nu=1}^rR_{\nu}$ follows from the description of non-wandering set of $f$. Take $R_{\nu}=N$. Since the boundary $\partial N$ consists of $(n-1)$-spheres, one can glue $n$-balls $B^n_1$, $\ldots$, $B^n_{j}$ to $\partial N$ to get a closed $n$-manifold $N_*$ where $j$ is the number of the boundary components of $N$. Due to Lemma \ref{lm:nbhd-of-attractor-omega}, the restriction $f|_{N}: N\to f(N)$ is extended to a Morse-Smale diffeomorphisms $f_*: N_*\to N_*$ whose non-wandering set consists of isolated nodal fixed points (sinks and sources). According to \cite{Reeb52} (see also \cite{GrinesMedvedevPochinkaZh2010}), $N_*$ is homeomorphic to $n$-sphere $\mathbb{S}^n$, and the non-wandering set of $f_*$ consists of a unique sink and unique source. Thus, $j=2$ and
$N=N_*\setminus (B^n_1\cup B^n_2)$. It follows from \cite{Keldysh1966} (see also \cite{DavermanVenema-book-2009}) that $N$ is homeomorphic to $\mathbb{S}^{n-1}\times [0;1]$. Moreover, one boundary component, say $\mathbb{S}^{n-1}\times\{0\}$, is a boundary component of some canonical neighborhood of expanding attractor while $\mathbb{S}^{n-1}\times\{1\}$ is a boundary component of some canonical neighborhood of contracting repeller.

For references, we formulate the following statement which can be obtained from \cite{MedvedevZhuzhoma2013-top-appl,MedvedevZhuzhoma2020-top-appl}.
For Reader's convenience, we give the sketch of proof.
\begin{proposition}\label{prop:ref-from-topol-appl}
Let $f: M^n\to M^n$ be a Morse-Smale diffeomorphism of closed $n$-manifold $M^n$, $n\geq 3$, such that $NW(f)$ consists of nodes and a unique saddle $\sigma$. If Morse index of $\sigma$ is different from $\{1,n-1\}$ then $NW(f)$ consists of a sink, a source, and the saddle $\sigma$. In addition, (1) $n\in\{4,8,16\}$; (2) $M^n$ is a simply connected manifold homeomorphic to $\mathbb{S}^{\frac{n}{2}}\sqcup\mathbb{B}^n$; (3) $M^n$ is a projective-like manifold provided $n\in\{8,16\}$.

If Morse index $ind(\sigma)$ of $\sigma$ belongs to $\{1,n-1\}$ then $M^n=\mathbb{S}^n$. Moreover, $NW(f)$ contains two sinks and a source provided $ind(\sigma)=1$ and $NW(f)$ contains a sink and two sources provided $ind(\sigma)=n-1$.
\end{proposition}
\textsl{Sketch of proof of Proposition \ref{prop:ref-from-topol-appl}}. If Morse index of $\sigma$ is different from $\{1,n-1\}$ then $f$ is a polar Morse-Smale diffeomorphism \cite{MedvedevZhuzhoma2020-top-appl}, Theorem 4. Therefore, $NW(f)$ consists of a sink, a source, and the saddle $\sigma$. It follows from \cite{MedvedevZhuzhoma2013-top-appl}, Theorem 1, that (1) $n\in\{4,8,16\}$; (2) $M^n$ is a simply connected manifold homeomorphic to $\mathbb{S}^{\frac{n}{2}}\sqcup\mathbb{B}^n$; (3) $M^n$ is a projective-like manifold provided $n\in\{8,16\}$.

It follows from \cite{OsenkovPochinka2024} for $n=3$ and \cite{MedvedevZhuzhoma2020-top-appl}, Theorem 1, for $n\geq 4$ that $M^n=\mathbb{S}^n$ provided $ind(\omega)\in\{1,n-1\}$. In addition, $NW(f)$ contains two sinks and a source provided $ind(\omega)=1$ and $NW(f)$ contains a sink and two sources provided $ind(\omega)=n-1$.
$\diamondsuit$

Suppose $f\in\mathbb{P}_k(M^n;0,0,1)$. Then the restriction $f|_{K_{\sigma}}: K_{\sigma}\to f(K_{\sigma})$ is extended to a Morse-Smale diffeomorphisms $\hat{f}: \widehat{K_{\sigma}}\to\widehat{K_{\sigma}}$ whose non-wandering set consists of isolated nodal fixed points and a unique saddle $\sigma$. Since $f\in\mathbb{P}_k(M^n;0,0,1)$, the Morse index of $\sigma$ is different from $\{1,n-1\}$. The result follows from Proposition \ref{prop:ref-from-topol-appl}.

Suppose $f\in\mathbb{M}_k(M^n;0,0,1)$. Similarly the restriction $f|_{K_{\sigma}}: K_{\sigma}\to f(K_{\sigma})$ is extended to a Morse-Smale diffeomorphisms $\hat{f}: \widehat{K_{\sigma}}\to\widehat{K_{\sigma}}$ whose non-wandering set consists of isolated nodal fixed points and a unique saddle $\sigma$. Since $f\in\mathbb{M}_k(M^n;0,0,1)$, the Morse index of $\sigma$ belongs to $\{1,n-1\}$. Again the result follows from Proposition \ref{prop:ref-from-topol-appl}.
$\Box$

\medskip
\textsl{Construction of graph $\Gamma_{\mathbb{P}}(f)$ for $f\in\mathbb{P}_k(M^n;0,0,1)$, $k\geq 2$, $n\in\{8,16\}$}.
Suppose $NW(f)$ consists of basic sets $\Omega_1$, $\ldots$, $\Omega_k$ with $l_1$, $\ldots$, $l_k$ bunches respectively. Due to Lemma \ref{lm:nbhd-of-attractor-omega}, there is a canonical neighborhood $U(\Omega_i)$ that is homeomorphic to $\mathbb{T}^n\setminus\cup_{j=1}^{l_i}D^n_j$, $i=1,\ldots,k$. Moreover, $f|_{U(\Omega_i)}$ can be extended to a DA-diffeomorphism $\mathbb{T}^n_i\to\mathbb{T}^n_i$ denoted by $f_i$ with a unique nontrivial basic set $\Omega_i$. Thus, $f_i$ corresponds the triple $(A_i,P_i,\epsilon_i)$. By definition, a graph $\Gamma_{\mathbb{P}}(f)$ has a group $V_i$ of vertices $v^i_1$, $\ldots$, $v^i_{l_i}$, and each $V_i$ endowed with the triple $(A_i,P_i,\epsilon_i)$ where $P_i$ consists of the points $p^i_1$, $\ldots$, $p^i_{l_i}$ (every $p^i_s$ corresponds $v^i_s$, and vice versa, $1\leq s\leq l_i$).

Now we keep the notation of Lemma \ref{lm:to-top-str-for-P-and-M}.
Consider a point $p^j_{t}\in P_j$, $t\in \{1, \dots k_j\}$, and a corresponding $(n-1)$-sphere $S^{n-1}(p^j_t)$ belonging to $\partial U(\Omega_j)$. We assume that the vertices $v^i_s$, $v^j_t$ are connected by an edge $L({v^i_s,v^j_t})$ provided the $(n-1)$-spheres $S^{n-1}(p^{i}_s)$, $S^{n-1}(p^j_t)$ bound a component $R_j$ with no non-wandering points. Sometimes, we'll say that the points $p^{i}_s$, $p^j_t$ are connected by the edge $L({p^i_s,p^j_t})$. It follows from Lemma \ref{lm:to-top-str-for-P-and-M} that $\epsilon_i\neq\epsilon_j$. If $S^{n-1}(p^{i}_s)$, $S^{n-1}(p^j_t)$ bound the component $K_{\sigma}$ then we assume that the vertices $v^i_s$, $v^j_t$ are connected by a \textit{marked edge} denoted by $e(f)$. This marked edge endowed with Pontryagin number $p^2_n(\widehat{K_{\sigma}})$ of the projective-like manifold $\widehat{K_{\sigma}}$ denoted by $p_e(f)$. Note that due to Kramer \cite{Kramer2003}, a topological type of $\widehat{K_{\sigma}}$ is completely determined by the Pontryagin number $p^2_n(\widehat{K_{\sigma}})$.

We see that $\Gamma_{\mathbb{P}}(f)$ is a collection of groups $V_1$, $\ldots$, $V_k$ of vertices such that the degree of every vertex equal to one, and there are no adjacent edges. Each group $V_i$ endowed with a triple $(A_i,P_i,\epsilon_i)$ corresponding to $\Omega_i$. Every edge $L({v^i_s,v^j_t})$ connects vertices with $\epsilon_i\neq\epsilon_j$, and there is a marked edge endowed with Pontryagin number. Roughly speaking, $\Gamma_{\mathbb{P}}(f)$ is a collection of pairwise disjoint segments. Below, we'll sometimes identify an automorphism with its matrix.

Suppose $\Gamma_{\mathbb{P}}(f)$, $\Gamma_{\mathbb{P}}(f')$ are graphs of diffeomorphisms $f,f'\in\mathbb{P}_k(M^n;0,0,1)$ respectively. We say that $\Gamma_{\mathbb{P}}(f)$, $\Gamma_{\mathbb{P}}(f')$ are \textit{commensurable} if the following conditions hold

(a) there is a bijection $\psi: \Gamma_{\mathbb{P}}(f)\to\Gamma_{\mathbb{P}}(f')$ such that $\psi(V_i)=V'_i$ for all $i=1,\ldots,k$, and $\psi(v^i_s)=v^{'i}_s$ for all $s=1,\ldots,l_i$. In particular, two vertices $v^i_s$, $v^j_t$ of $\Gamma(f)$ are connected by the edge $L({p^i_s,p^j_t})$ if and only if the vertices $\psi(v^i_s)=v^{'i}_s$, $\psi(v^j_t)=v^{'j}_t$ are connected by the edge $L({p^{'i}_s, p^{'j}_t})$. In addition, $\psi$ takes $p_e(f)$ to $p_e(f')$, and the marked edges $p_e(f)$, $p_e(f')$ have the same Pontryagin number;

(b) given any $1\leq i\leq k$, the triples $(A_i,P_i,\epsilon_i)$, $(A'_i,P'_i,\epsilon'_i)$ corresponding to the groups $V_i$, $V'_i$ are equivalent. Another words, there is a collection of automorphisms $\{\zeta_1,\ldots,\zeta_k\}$ of $\mathbb{T}^n$ such that $\zeta_i(P_i)=P'_i$, and $\zeta_i$ conjugates the automorphisms $A_i$, $A'_i$. In addition, the determinants of the automorphisms $\zeta_i$ are positive;

\medskip
\textsl{Construction of the set $\Gamma_{\mathbb{P}}^k$, $k\geq 2$, $n\in\{8,16\}$}.
Let us fix some $n\in\{8,16\}$. By definition, a graph $\gamma$ belongs to $\Gamma_{\mathbb{P}}^k$, $k\geq 2$, if it satisfies the following conditions :

(1) $\gamma\in\Gamma_{\mathbb{P}}^k$ has $k$ groups of vertices $V_i=\{v^i_1,\dots,v^i_{l_i}\}$, $i=1,\ldots,k$, and each group $V_i$ endowed with a triple $(A_i,P_i,\epsilon_i)$ where $A_i: \mathbb{T}^n\to\mathbb{T}^n$ is a codimension one Anosov automorphism, $P_i=\{p^i_1,\ldots,p^i_{l_i}\}$ is a finitely many fixed points of $A_i$, and $\epsilon_i=a$ provided the stable manifolds of $A_i$ is one-dimensional while $\epsilon_i=r$ provided the unstable manifolds of $A_i$ is one-dimensional. Moreover, there exists at least one triple $(A_i,P_i,\epsilon_i)$ with $\epsilon=a$ and at least one triple $(A_j,P_j,\epsilon_j)$ with $\epsilon=r$. In addition, there is a bijection $\psi_i: V_i\to P_i$, $p^i_s=\psi_i(v^i_s)$, for any $1\leq s\leq |P_i|=l_i$;

(2) every vertex $v^i_s$ of $\gamma\in\Gamma^k$ has degree 1, and there is a unique marked edge endowed with Pontryagin number $p^2_n$. For $n=8$, the Pontryagin number $p^2_4\in\{2(1+2t)^2 ; t\in\mathbb{Z}\}$, and for $n=16$, the Pontryagin number $p^2_8\in\{\frac{36}{49}(1+2t)^2 ; t\in\mathbb{Z}\}$;

(3) if vertices $v^i_s=\psi^{-1}_i(p^i_s)\in V_i=\{v^i_1,\dots,v^i_{l_i}\}$, $v^j_t=\psi^{-1}_j(p^j_t)\in V_j$ are connected by an edge $L({p^i_s,p^j_t})$, then $\epsilon_i\neq\epsilon_j$;

(4) given any groups $V_i$, $V_j$, $i\neq j$, there is a sequence of groups $V_{i_1}$, $l\dots$, $V_{i_r}$ with $i_1=i$, $i_r=j $ such that any neighbor groups $V_{i_s}$, $V_{i_{s+1}}$ in the sequence contain vertices
$v^{i_s}\in V_{i_s}$, $v^{i_{s+1}}\in V_{i_{s+1}}$ connected by an edge;

(5) the determinants of all $A_i$, $i=1,\ldots,k$, are positive.

\medskip
\textsl{Construction of graph $\Gamma_{\mathbb{M}}(f)$ for $f\in\mathbb{P}_k(M^n;0,0,1)$ and the set $\Gamma_{\mathbb{M}}^k$, $k\geq 2$.} Constructions of graphs $\Gamma_{\mathbb{M}}(f)$, $\Gamma_{\mathbb{P}}(f)$ differ just for marked edges. Due to Lemma \ref{lm:to-top-str-for-P-and-M}, $\widehat{K_{\sigma}}=K_{\sigma}\cup B_1^n\cup B_2^n\cup B_3^n$. Thus, $\Gamma_{\mathbb{M}}(f)$ has three marked vertices and two adjacent edges. Thus, one of the marked vertices has the degree two. Certainly, we omit the inclusion $n\in\{8,16\}$.

In the condition (a) of commensurability for $\mathbb{M}_k(M^n;0,0,1)$, we omit Pontryagin numbers for marked edges. In the condition (b), one requires that the determinants of the automorphisms $\zeta_i$ have the same sign.

As to the set $\Gamma_{\mathbb{M}}^k$, the condition (1) doesn't change. In condition (2), we require the existence of two adjacent marked edges with no Pontryagin numbers, and three marked vertices such that a unique marked point has the degree two. The conditions (3), (4) haven't changes. At last, in condition (5), one requires that the determinants of all $A_i$, $i=1,\ldots,k$, have the same sign.

\section{Topological structure of supporting manifolds}\label{s:topol}

In this section we revised the results of \cite{GrinMedvZh-2022}. We describe a topological structure of supporting manifolds for diffeomorphisms $f\in\mathbb{S}_k(M^n;a,b,c)\cup\mathbb{P}_k(M^n;a,b,c)\cup\mathbb{M}_k(M^n;a,b,c)$ because it is important for constructions of invariants of conjugacy.

The following result concerns to topological structure of supporting manifold $M^n$ for $f\in\mathbb{S}_k(M^n;a,b,c)$. Denote by $k_a$ (resp., $k_r$) the number of non-trivial attractors (resp., repellers) of $f$. We see that $k=k_a+k_r$.
\begin{lemma}\label{lm:top-str-for-S}
Let $f\in\mathbb{S}_k(M^n;a,b,c)$. Then $k+a+b=c+2$, and the supporting manifold $M^n$ is homeomorphic to
 $$ M^n=\underbrace{\mathbb{T}^n\sharp\cdots\sharp\mathbb{T}^n}_k. $$
Moreover, one of the following possibilities holds
\begin{itemize}
  \item either there is a unique codimension one contracting repeller ($k_r=1$), and there are no isolated sources ($b=0$), and $ind(f)=1$, $t_u(f)=1$, or
  \item there is a unique codimension one expanding attractor ($k_a=1$), and there are no isolated sinks ($a=0$), and $ind(f)=n-1$, $t_u(f)=k-1$.
\end{itemize}
\end{lemma}
\textsl{Proof}. Let $\Omega_1$, $\ldots$, $\Omega_k$ be nontrivial basic sets of $f$, and $U(\Omega_1)$, $\ldots$, $U(\Omega_k)$ its canonical neighborhoods respectively. Since $M^n$ is connected and every $\Omega_i$ is uniquely bunched, $K=M^n\setminus\left(\cup_{i=1}^kU(\Omega_i)\right)$ is a connected manifold with the boundary $\partial K=\cup_{i=1}^k\partial U(\Omega_i)$. Every $\partial U(\Omega_i)$ is an $(n-1)$-sphere. Therefore, one can glue to each component of $\partial K$ an $n$-ball $B^n_i$ to get a closed manifold $\widehat{K}=K\cup_{i=1}^kB^n_i$. Since any neighborhood $U(\Omega_i)$ is canonical, one can extend the restriction $f|_K$ to a diffeomorphism
$\hat{f}: \widehat{K}\to\widehat{K}$ which has in every ball $B^n_i$ a unique hyperbolic node, say $m_i$, $i=1,\ldots,k$. Moreover, $B^n_i\subset W^s(m_i)$ provided $m_i$ is a sink, and $B^n_i\subset W^u(m_i)$ provided $m_i$ is a source. Due to Condition ($2_{\mathbb{D}}$), $\hat{f}$ is a Morse-Smale diffeomorphism with $k+a+b$ nodes and $c$ saddles. It follows from Condition ($2_{\mathbb{S}}$) and \cite{MedvedevZhuzhoma2020-top-appl}, Theorem 3, that $\widehat{K}$ is an $n$-sphere. Hence, $M^n$ is a connected sum of $k$ copies of $\mathbb{T}^n$. Moreover, $k+a+b=c+2$.

Again, due to \cite{MedvedevZhuzhoma2020-top-appl}, Theorem 3, $k_r+b=1$ iff $ind(f)=1$ or $k_a+a=1$ iff $ind(f)=n-1$. Since $k_a\geq 1$ and $k_r\geq 1$, one gets $k_a=1$ or $k_r=1$. Clearly, $t_u(f)=1$ provided $k_r=1$, and $t_u(f)=k-1$ provided $k_a=1$.
$\Box$

\begin{cor}\label{cor:tube-agree}
Let $f\in\mathbb{S}_k(M^n;a,b,c)$. Then the $k$-tube $t(f)$ agreed with the triple $(a,b,c)$. To be precise, if $t(f)=(f_1,\ldots,f_k)$ contains a unique Anosov automorphism with one-dimensional unstable manifolds then $f$ has a unique codimension one contracting repeller and there are no isolated sources, and if $t(f)=(f_1,\ldots,f_k)$ contains a unique Anosov automorphism with one-dimensional stable manifolds then $f$ has a unique codimension one expanding attractor and there are no isolated sinks.
\end{cor}

\medskip
Topological structure of supporting manifolds for $f\in\mathbb{P}_k(M^n;0,0,1)$ is described as follows.
\begin{lemma}\label{top-struct-for-P}
Let $f\in\mathbb{P}_k(M^n;0,0,1)$, $k\geq 2$. Then the supporting manifold $M^n$ is homeomorphic to
 $$ M^n=\underbrace{\mathbb{T}^n\sharp\cdots\sharp\mathbb{T}^n}_{k\geq 2}\sharp\underbrace{\left({S}^{n-1}\times {S}^1\right)\sharp\cdots\sharp\left({S}^{n-1}\times {S}^1\right)}_{g\geq 0}\sharp N^n $$
where
\begin{itemize}
  \item $n\in\{4,8,16\}$
  \item $N^n$ is a simply connected manifold homeomorphic to $\mathbb{S}^{\frac{n}{2}}\sqcup\mathbb{B}^n$, here $\mathbb{B}^n$ is an open $n$-ball;
  \item $N^n$ is a projective-like manifold provided $n\in\{8,16\}$.
\end{itemize}
\end{lemma}
\textsl{Proof}. Let $\Omega_1$, $\ldots$, $\Omega_k$ be nontrivial basic sets of $f$, and $U(\Omega_1)$, $\ldots$, $U(\Omega_k)$ canonical neighborhoods respectively. It follows from Lemma \ref{lm:to-top-str-for-P-and-M} that $M^n\setminus\left(\cup_{i=1}^kU(\Omega_i)\right)$ is a disjoint union $K_{\sigma}\cup_{\nu=1}^rR_{\nu}$ where $K_{\sigma}$ contains $\sigma$ while $R_{\nu}$ is homeomorphic to $\mathbb{S}^{n-1}\times [0;1]$. Similarly to the proof of Lemma \ref{lm:top-str-for-S}, one can prove that $M^n$ is a connected sum of $k$ copies of $\mathbb{T}^n$, and copies of $\left({S}^{n-1}\times {S}^1\right)$, and $\widehat{K_{\sigma}}=K_{\sigma}\cup B_1^n\cup B_2^n$. The result follows immediately from Lemma \ref{lm:to-top-str-for-P-and-M}.
$\Box$

\begin{lemma}\label{top-struct-for-M}
Let $f\in\mathbb{M}_k(M^n;0,0,1)$, $k\geq 2$. Then the supporting manifold $M^n$ is homeomorphic to
 $$ M^n=\underbrace{\mathbb{T}^n\sharp\cdots\sharp\mathbb{T}^n}_{k\geq 2}\sharp\underbrace{\left({S}^{n-1}\times {S}^1\right)\sharp\cdots\sharp\left({S}^{n-1}\times {S}^1\right)}_{g\geq 0}. $$
\end{lemma}
\textsl{Proof}. We keep the notation of Lemmas \ref{lm:to-top-str-for-P-and-M}, \ref{top-struct-for-P}. Here, $\widehat{K_{\sigma}}=K_{\sigma}\cup B_1^n\cup B_2^n\cup B_3^n=\mathbb{S}^n$. As a consequence,$M^n$ is a connected sum of $k$ copies of $\mathbb{T}^n$, and copies of $\left({S}^{n-1}\times {S}^1\right)$.
$\Box$

\section{Proofs of the main results}\label{s:proofs}

\textsl{Proof of Theorem \ref{thm:classif-for-S}}. \textbf{Necessity.} Suppose diffeomorphisms $f_i\in\mathbb{S}_k(M^n;a_i,b_i,c_i)$, $i=1,2$, are globally conjugate on their non-wandering sets, i.e. there is a homeomorphism
$\varphi: M^n\to M^n$ such that $\varphi\left(NW(f_1)\right)=NW(f_2)$ and $\varphi\circ f_1|_{NW(f_1)}=f_2\circ\varphi|_{NW(f_1)}$. Then $\varphi$ induces a one-to-one correspondence between the non-trivial basic sets $\Omega^1_1$, $\ldots$, $\Omega^1_k$ of $f_1$ and nontrivial basic sets $\Omega^2_1$, $\ldots$, $\Omega^2_k$ of $f_2$, i.e. there is a $k$-permutation $\tau$ such that $\varphi(\Omega^1_i)=\Omega^2_{\tau(i)}$, $i=1,\ldots,k$. For simplicity, we'll assume $\tau(i)=i$, $i=1,\ldots,k$.

Take pairwise disjoint canonical neighborhoods $U(\Omega^j_1)$, $\ldots$, $U(\Omega^j_k)$, $j=1,2$. Without loss of generality, one can assume that $\varphi$ takes any canonical neighborhood $U(\Omega^1_i)$ to canonical neighborhood $U(\Omega^2_i)$. Hence, $\varphi$ induced the conjugacy of corresponding DA-diffeomorphisms. According Theorem \ref{thm:plyk-gr-class}, for every $f_j$, one gets the $k$-tube $t(f_j)=(f^j_1,\ldots,f^j_k)$, $j=1,2$. As a consequence, the $k$-tubes $t(f_1)$, $t(f_2)$ are equivalent under some linear automorphisms $\zeta_1$, $\ldots$, $\zeta_k$. later on, we'll need the following result.
\begin{proposition}\label{prop:sphere-with-holes-homeo}
Let $B^n_1$, $\ldots$, $B^n_k$ be pairwise disjoint open $n$-balls in $\mathbb{S}^n$, so that $K^n=\mathbb{S}^n\setminus\left(\cup_{i=1}^kB^n_i\right)$ is a connected manifold with the boundary
$\partial K^n=\cup_{i=1}^k\partial B^n_i$. Suppose a homeomorphism $\varphi: K^n\to K^n$ takes every $\partial B^n_i$ to itself, $i=1,\ldots,k$. Then $\varphi$ preserves orientation iff every restriction
$\varphi|_{\partial B^n_i}: \partial B^n_i\to\partial B^n_i$ preserves orientation (so, $\varphi$ reverses orientation iff every restriction $\varphi|_{\partial B^n_i}$ reveres orientation), $i=1,\ldots,k$. Conversely, let $\phi_i: \partial B^n_i\to\partial B^n_i$ be a homeomorphism, $i=1,\ldots,k$. Then there is a homeomorphism $\varphi: K^n\to K^n$ such that $\varphi|_{\partial B^n_i}$, $i=1,\ldots,k$, iff one of the following possibilities holds:
1) all $\phi_i$, $i=1,\ldots,k$, preserve orientation; 2) all $\phi_i$, $i=1,\ldots,k$, reverse orientation.
\end{proposition}
\textit{Proof of Proposition \ref{prop:sphere-with-holes-homeo}}. For $k=1$, the manifold $K^n$ is a closed $n$-disk, and the result is obvious \cite{Hatcher-book}, Section 3.3. For $k=2$, the manifold $K^n$ is a closed $n$-annulus $\mathbb{S}^{n-1}\times [0;1]$. Homeomorphisms $\phi_1$, $\phi_2$ are extended to $K^n$ iff they are isotopic. According to \cite{GrinGurevich-book-2022}, Theorem 14.5, $\phi_1$, $\phi_2$ either the both preserve orientation or the both reverse orientation. Now the proof is by induction over $k$. We see that for $k=1,2$, the result is true. Suppose that the result is true for $1,2,\ldots,k$, and let us prove the statement for $k+1$. If $\varphi: K^n\to K^n$ is a homeomorphisms then the both $\varphi$ and $\varphi|_{\partial B^n_{k+1}}$ either preserve or reverse orientation since $\partial B^n_{k+1}$ bounds the ball $B^n_{k+1}$ \cite{Hatcher-book}, Section 3.3. Let now $\phi_i: \partial B^n_i\to\partial B^n_i$ be a homeomorphism, $i=1,\ldots,k,k+1$. By the induction hypothesis, the homeomorphisms $\phi_1$, $\ldots$, $\phi_k$ can be extended to a homeomorphism $\varphi: K^n\cup B^n_{k+1}\to K^n\cup B^n_{k+1}$ iff one of the possibilities above holds. Since the balls $B^n_i$ be a homeomorphism, $i=1,\ldots,k$, $\varphi(B^n_{k+1})$ are pairwise disjoint, one can deform $\varphi$ such that $\varphi(B^n_{k+1})=B^n_{k+1}$. Suppose all $\phi_i$, $i=1,\ldots,k,k+1$, and $\varphi|_{\partial B^n_{k+1}}$ preserve orientation. Then $\phi_{k+1}: \partial B^n_{k+1}\to\partial B^n_{k+1}$ and $\varphi|_{\partial B^n_{k+1}}$ are isotopic \cite{GrinGurevich-book-2022}, Theorem 14.5. Therefore one can deform $\varphi$ in a neighborhood of $\partial B^n_{k+1}$ to get $\varphi$ such that $\varphi|_{\partial B^n_{k+1}}=\phi_{k+1}$. Similarly, one considers the case when all $\phi_i$, $i=1,\ldots,k,k+1$, and $\varphi|_{\partial B^n_{k+1}}$ reverse orientation. This completes the proof.
$\diamondsuit$

\medskip
It follows from Lemma \ref{lm:orient-determinant} and Proposition \ref{prop:sphere-with-holes-homeo} that the determinants of all $\zeta_j$ have the same sign.

Obviously, $a_1+b_1+c_1=a_2+b_2+c_2$. It follows from Lemma \ref{lm:top-str-for-S} that $k+a_1+b_1=c_1+2$ and $k+a_2+b_2=c_2+2$. Hence, $c_1=c_2$ and $a_1+b_1=a_2+b_2$.

\medskip
\textbf{Sufficiency.} Suppose the $k$-tubes $t(f_1)=(A_1^1,\ldots,A_k^1)$, $t(f_2)=(A_1^2,\ldots,A_k^2)$ of $f_1$, $f_2$ respectively are equivalent under linear automorphisms $\zeta_1$, $\ldots$, $\zeta_k$, and the determinants of all $\zeta_j$ have the same sign. Without loss of generality, one can assume that $\zeta_j$ conjugates $A_j^1$, $A_j^2$, $j=1,\ldots,k$. Since $f_i\in\mathbb{S}_k(M^n;a_i,b_i,c_i)$, the corresponding DA-diffeomorphisms $(f_1^i,\ldots,f_k^i)$, $i=1,2$, are simplest ones. It follows from Lemma \ref{lm:conjugacy-of-simplest-da} that $f_j^1$, $f_j^2$ are globally conjugate on the basic sets $\Omega_j^1$, $\Omega_j^2$, $j=1,\ldots,k$, i.e. there is a homeomorphism $\phi_j: \mathbb{T}^n_j\to\mathbb{T}^n_j$ such that $\phi_j(\Omega_j^1)=\Omega_j^2$ and $f_j^1\circ\phi_j=\phi_j\circ f_j^2$ for every $j$. Here, $\mathbb{T}^n_j$ is a copy of $\mathbb{T}^n$. Since the determinants of all $\zeta_j$ have the same sign, all $\phi_j$ can be extended to a common homeomorphism $\varphi: M^n\to M^n$ due to Proposition \ref{prop:sphere-with-holes-homeo}. Since $a_1+b_1=a_2+b_2$, $c_1=c_2$, one can assume that $\varphi$ takes all isolated nodes and saddles of $f_1$ to isolated nodes and saddles of $f_2$ respectively.

Suppose $f\in\mathbb{S}_k(M^n;a,b,c)$. According to Lemma \ref{lm:top-str-for-S}, $k+a+b=c+2$ and $M^n$ of the kind (\ref{eq:S-decomposition-manifold}). It follows from Lemmas \ref{lm:nbhd-of-attractor-omega}, \ref{lm:top-str-for-S} that the $k$-tube $k(f)$ is admissible, $t(f)\in A_k$.

\medskip
\textbf{Realization.}
Take $k$-tube $t=(F_1,\ldots,F_k)\in A_k$ and integers $k\geq 2$, $a,b\geq 0$, $c\geq 1$ such that $k+a+b=c+2$. Assume for definiteness that there is exactly one Anosov automorphism, say $F_1$, in $t$ with a stable one-dimension manifolds. According to Franks \cite{Fr69}, every $F_i$ has a fixed point. Therefore one can construct simplest DA-diffeomorphisms $f_1$, $\ldots$, $f_k$ such that $f_{i*}=F_i$, $i=1,\ldots,k$ where $f_1$ has a codimension one expanding attractor while another $f_j$, $j=2,\ldots,k$, have a codimension one contracting repeller \cite{GinesZhuzhoma2005,Plykin84}. Since $t$ agreed with the triple $(a,b,c)$ and $t_u=k-1$, $a=0$. It follows from Lemma \ref{lm:top-str-for-S} that isolated saddles must have the Morse index $ind(f)=n-1$.

There is a Morse-Smale diffeomorphism $g: \mathbb{S}^n\to\mathbb{S}^n$ with $NW(g)$ consisting of a unique sink, say $\omega_0$, $k+b$ sources, and $c=k+b-2$ saddles with the Morse index $n-1$ \cite{MedvedevZhuzhoma2020-top-appl}, Section 3. The diffeomorphism $f_1$ has an expanding attractor, say $\Lambda_a$, and a unique source, say $\alpha_0$. Deleting neighborhoods of $\omega_0$ and $\alpha_0$, one can construct on the manifold $\mathbb{T}^n\sharp\mathbb{S}^n=\mathbb{T}^n$ a diffeomorphism with the non-wandering set $\left(NW(g)\setminus\{\omega_0\}\right)\cup\Lambda_a$. Similarly one can replace $k-1$ isolated sources by $k-1$ nontrivial repellers.
This completes the proof.
$\Box$

\medskip
\textsl{Proof of Theorem \ref{thm:classif-for-P}}.
\textbf{Necessity.} Suppose diffeomorphisms $f_1,f_2\in\mathbb{P}_k(M^n;0,0,1)$, $k\geq 2$ and $n\in\{8,16\}$, are globally conjugate on their non-wandering sets, i.e. there is a homeomorphism
$\varphi: M^n\to M^n$ such that $\varphi\left(NW(f_1)\right)=NW(f_2)$ and $\varphi\circ f_1|_{NW(f_1)}=f_2\circ\varphi|_{NW(f_1)}$. The conjugacy on nontrivial basic sets implies that $\varphi$ takes the boundary points of $f_1$ onto the boundary points of $f_2$. Hence, $\varphi$ induces a one-to-one correspondence between the non-trivial basic sets $\Omega^1_1$, $\ldots$, $\Omega^1_k$ of $f_1$ and nontrivial basic sets $\Omega^2_1$, $\ldots$, $\Omega^2_k$ of $f_2$, i.e. there is a $k$-permutation $\tau$ such that $\varphi(\Omega^1_i)=\Omega^2_{\tau(i)}$, $i=1,\ldots,k$. For simplicity, we'll assume $\tau(i)=i$, $i=1,\ldots,k$. As a consequence, $\varphi$ induces a bijection $\psi$ of groups $V_i=\{v^i_1,\dots,v^i_{k_i}\}\to V_i'=\{v^{'i}_1,\dots,v^{'i}_{k_i}\}$, $i=1,\ldots,k$, and vertices $\psi(v^i_s)=v^{'i}_s$, $s=1,\ldots,k_i$ inside of groups. Due to Lemma \ref{lm:to-top-str-for-P-and-M}, $\psi$ takes $p_e(f_1)$ to $p_e(f_2)$ keeping its Pontryagin number. We see that the condition (а) of commensurability  holds.

Take pairwise disjoint canonical neighborhoods $U(\Omega^j_1)$, $\ldots$, $U(\Omega^j_k)$, $j=1,2$. Without loss of generality, one can assume that $\varphi$ takes any canonical neighborhood $U(\Omega^1_i)$ to canonical neighborhood $U(\Omega^2_i)$. According Lemma \ref{lm:nbhd-of-attractor-omega}, the restriction $f_j|_{U(\Omega^j_i)}$ is extended to DA-diffeomorphism $\tilde{f^j_i}: \mathbb{T}^n\to\mathbb{T}^n$, $j=1,2$, $i=1,\ldots,k$.
Therefore, $\varphi$ induced the conjugacy of the DA-diffeomorphisms $\tilde{f^1_i}$, $\tilde{f^2_i}$. We see that $\tilde{f^j_i}$ semi-conjugates to the codimension one Anosov automorphism
$(\tilde{f^j_i})_*: \mathbb{T}^n\to\mathbb{T}^n$ denoted by $f^j_i$. As a consequence, the triples $(f^1_i,P^1_i,\epsilon^1_i)$, $(f^2_i,P^2_i,\epsilon^2_i)$ are equivalent. Thus, there are linear automorphisms $\zeta_1$, $\ldots$, $\zeta_k$ such that $f^1_i$ conjugates $f^2_i$ under $\zeta_i$, $i=1,\ldots,k$.

Now we'll keep the notation of Lemma \ref{lm:to-top-str-for-P-and-M}. Since $f_1,f_2\in\mathbb{P}_k(M^n;0,0,1)$ and $n\in\{8,16\}$, the manifolds $\widehat{K_{\sigma_j}}$, $j=1,2$, are projective-like. Clearly, the restriction $\varphi|_{K_{\sigma_1}}: K_{\sigma_1}\to K_{\sigma_2}$ is extended to a homeomorphism $\widehat{K_{\sigma_1}}\to\widehat{K_{\sigma_2}}$ denoted by $\phi_*$. Due to \cite{EellsKuiper62,Kramer2003}, Pontryagin numbers of $\widehat{K_{\sigma_j}}$, $j=1,2$, are nonzero. It follows from \cite{Novikov-1965} that $\phi_*$ is an orientation preserving mapping. This implies that the restriction of $\varphi$ on the boundary components of $K_{\sigma_1}$ preserve orientation also. According to Lemma \ref{lm:orient-determinant}, the determinants of all $\zeta_i$ are positive since the manifold $M^n$ is connected. We see that the condition (b) of commensurability  holds. Thus, the graphs $\Gamma_{\mathbb{P}}(f_1)$, $\Gamma_{\mathbb{P}}(f_2)$ are commensurable.

\textbf{Sufficiency.} Suppose the graphs $\Gamma_{\mathbb{P}}(f_1)$, $\Gamma_{\mathbb{P}}(f_2)$ are commensurable. Due to the condition (a), there is a bijection $\psi: \Gamma_{\mathbb{P}}(f_1)\to\Gamma_{\mathbb{P}}(f_2)$ which induces a bijection of groups of vertices. Recall that every group of vertices corresponds a unique basic set. Therefore, $\psi$ induces a one-to-one bijection $\Omega_i\Longleftrightarrow\Omega'_i$, $i=1,\ldots,k$, between the basic sets of $f_1$, $f_2$. According the condition (b), given any $1\leq i\leq k$, the triples $(A_i,P_i,\epsilon_i)$, $(A'_i,P'_i,\epsilon'_i)$ correspon\-ding to the groups $V_i=\{v^i_1,\ldots,v^i_{l_i}\}$, $V'_i=\{v^{'i}_1,\ldots,v^{'i}_{l_i}\}$ are equivalent. Another words, there is a collection of automorphisms $\{\zeta_1,\ldots,\zeta_k\}$ of $\mathbb{T}^n$ such that $\zeta_i(P_i)=P'_i$, and $\zeta_i$ conjugates the automorphisms $A_i$, $A'_i$, $i=1,\ldots,k$. Therefore there is a homeomorphism $\varphi_i: U(\Omega_i)\to U(\Omega_i')$ such that $\varphi_i\circ f|_{\Omega_i}=f'\circ \varphi_i|_{\Omega_i}$. We have to prove that the homeomorphisms $\varphi_i: U(\Omega_i)\to U(\Omega_i')$, $i=1,\ldots,k$, can be extended to a common homeomorphism $\varphi: M^n\to M^n$.

Without loss of generality, one can assume that $\psi(v^i_s) = v^{'i}_s$ for all $i=1,\ldots,k$ and $s=1,\ldots,l_i$ where $V_i=\{v^i_1,\dots,v^i_{l_i}\}$ are groups of the vertices of $\Gamma_{\mathbb{P}}(f_1)$, and $V_i'=\{v^{'i}_1,\dots,v^{'i}_{l_i}\}$ are groups of vertices of $\Gamma_{\mathbb{P}}(f_2)$, $i=1,\ldots,k$. Suppose non-marked vertices $v^i_s$, $v^j_t\in\Gamma_{\mathbb{P}}(f_1)$ is connected by an edge $L({p^i_s, p^j_t})$. Hence the vertices $\psi(v^i_s)=v^{'i}_s$, $\psi(v^j_t)=v^{'j}_t\in\Gamma_{\mathbb{P}}(f_2)$ is connected by the non-marked edge $L({p^{'i}_s, p^{'j}_t})$ also.
By Lemma \ref{lm:to-top-str-for-P-and-M}, $(n-1)$-spheres $S^{n-1}(p^{i}_s)$, $S^{n-1}(p^j_t)$ bound a set $K_{ij}$ homeomorphic to $S^{n-1}\times [0;1]$. Similarly, the $(n-1)$-spheres $\varphi_i(S^{n-1}(p^{i}_s))=S^{n-1}(p^{'i}_s)$, $\varphi_j(S^{n-1}(p^{j}_t)=S^{n-1}(p^{'j}_t)$ bound an $n$-annulus $K_{ij}'\subset M^n$ homeomorphic to $S^{n-1}\times [0;1]$. Since the automorphisms $\zeta_i$, $\zeta_j$ have the same sign, $\varphi_i$ and $\varphi_j$ can be extended to a homeomorphism $\varphi_{ij}: K_{ij}\to K_{ij}'$.

Due to the condition (a), $\psi$ takes $p_e(f)$ to $p_e(f')$, and the marked edges $p_e(f)$, $p_e(f')$ have the same Pontryagin number. This follows that the projective-like manifolds $\widehat{K_{\sigma_1}}$, $\widehat{K_{\sigma_2}}$ are homeomorphic \cite{Kramer2003}. Since the determinants of the corresponding automorphisms $\zeta_i$, $\zeta_j$ are positive, one can extend $\varphi_i$ and $\varphi_j$ to a homeomorphism $\widehat{K_{\sigma_1}}\to\widehat{K_{\sigma_2}}$. Clearly, such homeomorphism can be deformed to a homeomorphism which takes $\sigma_1$ to $\sigma_2$. We see that the homeomorphisms $\varphi_i$ can be extended to a common homeomorphism $\varphi: M^n\to M^n$. Thus, $f$ and $f'$ are globally conjugate on its non-wandering sets.

\medskip
\textbf{Realization.} First, let us show that a graph $\Gamma_{\mathbb{P}}(f)$ belongs to $\Gamma^k_{\mathbb{P}}$. It follows from the construction of $\Gamma_{\mathbb{P}}(f)$ that given any basic set $\Omega_i$ with $l_i\geq 1$ bunches, one corresponds a group $V_i=\{v^i_1,\dots,v^i_{l_i}\}$ of vertices endowed with a triple $(A_i,P_i,\epsilon_i)$ where $A_i: \mathbb{T}^n\to\mathbb{T}^n$ is a codimension one Anosov automorphism, $P_i=\{p^i_1,\ldots,p^i_{l_i}\}$ is an invariant set of $A_i$ consisting of finitely many periodic points. Moreover, $\epsilon_i=a$ provided the stable manifolds of $A_i$ is one-dimensional while $\epsilon_i=r$ provided the unstable manifolds of $A_i$ is one-dimensional. It follows from Lemma \ref{lm:to-top-str-for-P-and-M} that there exists at least one triple $(A_i,P_i,\epsilon_i)$ with $\epsilon=a$ and at least one triple $(A_j,P_j,\epsilon_j)$ with $\epsilon=r$. Since $|P_i|=l_i$, there is a one-to-one correspondence between points of $P_i$ and vertices of $V_i$. Without loss of generality, one can assume that there is a bijection $\psi_i: V_i\to P_i$ such that $p^i_s=\psi_i(v^i_s)$, $i=1,\ldots,k$, $s=1,\ldots,k_i$. Moreover, there is a marked edge endowed with Pontryagin number $p^2_n$. Due to \cite{Kramer2003}, for $n=8$, the Pontryagin number $p^2_4\in\{2(1+2t)^2 ; t\in\mathbb{Z}\}$, and for $n=16$, the Pontryagin number $p^2_8\in\{\frac{36}{49}(1+2t)^2 ; t\in\mathbb{Z}\}$. We see that the conditions (1) and (2) of the description of the set $\Gamma^k$ hold. Again, Lemma \ref{lm:to-top-str-for-P-and-M} implies the condition (3). Since $M^n$ is a connected manifold, the condition (4) holds as well.

Suppose vertices $v^i_s$, $v^j_t\in\Gamma(f)$ are connected by a non-marked edge $L({p^i_s,p^j_t})$. Recall that the vertices $v^i_s$, $v^j_t$ belong to some groups endowed with the triples $(A_i,P_i,\epsilon_i)$, $(A_j,P_j,\epsilon_j)$ respectively, so that $p^{i}_s\in P_i$, $p^j_t\in P_j$. Since the edge $L({p^i_s,p^j_t})$ is non-marked, the spheres $S^{n-1}(p^i_s)$, $S^{n-1}(p^j_t)$ bound an $n$-annulus $K_{ij}\subset M^n$ where $S^{n-1}(p^i_s)$, $S^{n-1}(p^j_t)$ are components of the boundary of attracting neighborhoods $U(\Omega_i)$, $U(\Omega_j)$ of basic sets $\Omega_i$, $\Omega_j$ respectively. Since $S^{n-1}(p^i_s)$, $S^{n-1}(p^j_t)$ bound the annulus $K_{ij}$ that is homeomorphic to $S^{n-1}\times [0;1]$, the restrictions $f|_{S^{n-1}(p^i_s)}$, $f|_{S^{n-1}(p^j_t)}$ are isotopic. Due to Proposition \ref{prop:sphere-with-holes-homeo}, this restrictions $f|_{S^{n-1}(p^i_s)}$, $f|_{S^{n-1}(p^j_t)}$ are either the both preserve orientation or the both reverse orientation. According to Lemma \ref{lm:nbhd-of-attractor-omega}, the diffeomorphisms $f_i$, $f_j$ are extended to DA-diffeomorphisms $\tilde{f}_i$, $\tilde{f}_j: \mathbb{T}^n\to\mathbb{T}^n$ respectively. Moreover, the both $S^{n-1}(p^i_s)$ and $S^{n-1}(p^j_t)$ bounds an $n$-ball in $\mathbb{T}^n$. Since $f|_{S^{n-1}(p^i_s)}$, $f|_{S^{n-1}(p^j_t)}$ are isotopic, the diffeomorphisms $\tilde{f}_i$, $\tilde{f}_j$ are either the both preserve orientation or the both reverse orientation. It follows from the relations $A_i=(\tilde{f}_i)_*$, $A_j=(\tilde{f}_j)_*$ that the determinants of $A_i$, $A_j$ have the same sign. Since $f\in\in\mathbb{P}_k(M^n;0,0,1)$ and $n\in\{8,16\}$, there is a component of $M^n\setminus\left(\cup_{i=1}^kU(\Omega_i)\right)$ such that $\widehat{K_{\sigma}}=K_{\sigma}\cup B_1^n\cup B_2^n$ is a projective-like $n$-manifold. Due to \cite{Novikov-1965} and Kramer \cite{Kramer2003}, a projective-like manifold admit only preserving orientation mapping. Hence, $f|_{\partial K_{\sigma}}$ preserve orientation. This implies that corresponding $A_i$, $A_j$ preserve orientation. As a consequence, the determinants of $A_i$, $A_j$ are positive. Now, the condition (5) follows from (4).

Take now an abstract graph $\gamma\in\Gamma^k_{\mathbb{P}}$. We have to show that there are a closed orientable $n$-manifold $M^n$, $n\in\{8,16\}$, and a diffeomorphism $f\in\mathbb{P}_k(M^n;0,0,1)$, $k\geq 2$, such that $\gamma=\Gamma_{\mathbb{P}}(f)$. It follows from condition (1) that every group of vertices $V_i=\{v^i_1,\dots,v^i_{k_i}\}$ endowed with a triple $(A_i,P_i,\epsilon_i)$ where $A_i: \mathbb{T}^n\to\mathbb{T}^n$ is Anosov automorphism with a finite invariant set of fixed points $P_i$, and $\epsilon_i=a$ provided the stable manifolds of $A_i$ is one-dimensional while $\epsilon_i=r$ provided the unstable manifolds of $A_i$ is one-dimensional.
Using Smale's surgery operation \cite{Smale67} (see also \cite{GinesZhuzhoma2005,Plykin84}), one can construct a DA-diffeomorphism $f_i: \mathbb{T}^n\to\mathbb{T}^n$ with codimension one orientable connected basic set $\Omega_i$ containing $|P_i|=k_i$ bunches. Moreover, if $\epsilon_i= a$ then $\Omega_i$ is an expanding attractor, and if $\epsilon_i=r$ then $\Omega_i$ is a contracting repeller. In addition, every bunch corresponds some point
$p^i_s\in P_i$ and vertex $v^i_s=\psi^{-1}(p^i_s)\in V_i$. According \cite{GinesZhuzhoma1979,Plykin84,GinesZhuzhoma2005}, the triple $(A_i,P_i,\epsilon_i)$ is a complete invariant of conjugacy for the diffeomorphism $f_i$. Recall that every component of $\mathbb{T}^n\setminus\Omega_i$ contains a unique isolated node fixed point surrounded by a characteristic sphere of corresponding bunch. Hence, all boundary points of $f_i$ are fixed.

Let us take $k$ copies $\mathbb{T}^n_1$, $\ldots$, $\mathbb{T}_k^n$ of $\mathbb{T}^n$. It is convenient to consider $f_i: \mathbb{T}^n_i\to\mathbb{T}^n_i$ defined on $\mathbb{T}^n_i$, $i=1,\ldots,k$. Due to condition (2), every vertex $v^i_s=\psi^{-1}(p^i_s)$ is connected with a unique vertex $v^j_t=\psi^{-1}_j(p^j_t)$ by an edge $L(p^i_s,p^j_t)\subset\gamma$ where $i\neq j$. According to condition (3), the vertex $v^j_t$ belongs to a group $V_j$ endowed with a triple $(A_j,P_j,\epsilon_j)$ where $p^j_t=\psi_j(v^j_t)\in P_j$, $\epsilon_i\neq \epsilon_j$. For definiteness, assume that $\epsilon_i=a$ and $\epsilon_j=r$, i.e. $\Omega_i$ is an attractor and $\Omega_j$ is a repeller.

Due to \cite{Plykin84}, there is an attracting neighborhood $U(\Omega_i)\subset W^s(\Omega_i)$ of $\Omega_i$ such that the set $\mathbb{T}_i^n\setminus U(\Omega_i)$ is the union of pairwise disjoint $n$-disks $B^i_1$, $\ldots$, $B^i_{k_i}$, and the boundary $\partial U(\Omega_i)$ is the union of characteristic spheres $\widehat{S}^i_1=\partial B^i_1$, $\ldots$, $\widehat{S}^i_{k_i}=\partial B^i_{k_i}$, $s=1,\ldots,k_i$. Since $f_i$ is a DA-diffeomorphism, every $n$-disk $B^i_m$ contains a unique source fixed point. Some source fixed point denoted by $q^i_s$ corresponds to the vertex $v^i_s$. Without loss of generality, one can assume that $q^i_s\in B^i_1$. Since $f_i\left(U(\Omega_i)\right)\subset U(\Omega_i)$, $B^i_1\subset f(B^i_1)$. Similarly, there is an attracting neighborhood $U(\Omega_j)\subset W^u(\Omega_j)$ of the basic set $\Omega_j$ such that the set
$\mathbb{T}_j^n\setminus U(\Omega_j)$ is the union of pairwise disjoint $n$-disks $B^j_1$, $\ldots$, $B^j_{k_j}$, and the boundary $\partial U(\Omega_j)$ is the union of characteristic spheres $\widehat{S}^j_1$, $\ldots$, $\widehat{S}^j_{k_j}$. Every $n$-disk $B^j_m$ contains a unique sink fixed point. sink fixed point denoted by $q^j_t$ corresponds to the vertex $v^j_t$. Without loss of generality, one can assume that $q^j_t\in B^j_1$.

Suppose the edge $L(p^i_s,p^j_t)$ is non-marked. Let us delete the disks $B^i_1$, $f_j(B^j_1)$ from $\mathbb{T}^n_i$, $\mathbb{T}^n_j$ respectively. Take an $n$-annulus $K_{ij}$, and glue its boundary component to
$\partial B^i_1$, $\partial f_j(B^j_1)$ so that the set
 $$\left(\mathbb{T}^n_i\setminus B^i_1\right)\bigcup\left(\mathbb{T}^n_j\setminus f_j(B^j_1)\right)\bigcup K_{ij}=M^n_{ij}$$
becomes a smooth closed orientable manifold. Since $B^i_1\subset f_i(B^i_1)$ and $f_j(B^j_1)\subset B^j_1$, the topological closures $K_i$, $K_j$ of the sets $f_i(B^i_1)\setminus B^i_1$, $B^j_1\setminus f_j(B^j_1)$ respectively are $n$-annuluses. Therefore, $K_{ij}\cup K_i$ is an $n$-annulus with two boundary components $\partial f_i(B^i_1)$, $\partial f_j(B^j_1)$ while the union $K_{ij}\cup K_j$ is an $n$-annulus with the boundary components
$\partial B^i_1$, $\partial B^j_1$. By condition (5), the determinants of $A_i$, $A_j$ have the same sign. This implies that the restrictions $f_i|_{\partial B^i_1}$, $f_j|_{\partial B^j_1}$ either the both preserve orientation or the both reverse orientation. Due to \cite{GrinGurevich-book-2022}, Theorem 14.5, this restrictions are isotopic. Therefore, there is a mapping
 $$ \varphi_{ij}: K_{ij}\cup K_j\to K_{ij}\cup K_i \, \mbox{ such that }\, \varphi_{ij}|_{\partial B^i_1}=f_i|_{\partial B^i_1},\,\, \varphi_{ij}|_{\partial B^j_1}=f_j|_{\partial B^j_1}. $$
Since $K_{ij}\cup K_j$ is an $n$-annulus, we can define $\varphi_{ij}$ so that all points on $K_{ij}\cup K_j\cup K_i$ move from $\partial B^j_1$ to $\partial B^i_1$ under positive iteration of $\varphi_{ij}$. Moreover, one can $\varphi_{ij}$ agrees with the restrictions $f_i|_{\partial B^i_1}$, $f_j|_{\partial B^j_1}$ near $\partial B^i_1$, $\partial B^j_1$ respectively so that a mapping
 $$ f_{ij}|_{\mathbb{T}^n_i\setminus B^i_1}=f_i|_{\mathbb{T}^n_i\setminus B^i_1},\,\, f_{ij}|_{\mathbb{T}^n_j\setminus B^j_1}=f_j|_{\mathbb{T}^n_j\setminus B^j_1},\,\, f_{ij}|_{K_{ij}\cup K_j}=\varphi_{ij} $$
becomes a diffeomorphism $f_{ij}: M^n_{ij}\to M^n_{ij}$. Taking in mind the property of $\varphi_{ij}|_{K_{ij}\cup K_j}$, we see that $f_{ij}$ has no non-wandering points on $K_{ij}\cup K_j\cup K_i$. Hence, $f_{ij}$ is an A-diffeomorphism whose non-wandering set consists of the orientable connected codimension basic sets $\Omega_i$, $\Omega_j$ and trivial basic sets of the diffeomorphisms $f_i$, $f_j$ without the points $q^i_s$, $q^j_t$. The $n$-annulus $K_{ij}$, in sense, corresponds to the edge $L(p^i_s,p^j_t)$.

Suppose now that the edge $L(p^i_s,p^j_t)$ is marked where $v^i_s=\psi^{-1}(p^i_s)\in V_i$ and $v^i_s=\psi^{-1}(p^j_t)\in V_j$. Denote by $S^{n-1}(p^i_s)$ (resp., $S^{n-1}(p^j_t)$) the boundary component of $U(\Omega_i)$
(resp., $U(\Omega_j)$) corresponding to the point $p^i_s$ (resp., $p^j_t$). For definiteness, assume that $\Omega_i$ is an attractor while $\Omega_j$ is a repeller.
By condition (2), $L(p^i_s,p^j_t)$ endowed with Pontryagin number $p^2_n$. According \cite{EellsKuiper62,Kramer2003}, there is a projective-like manifold $N^n$ such that the Pontryagin number of $N^n$ equals $p^2_n$. Due to \cite{MedvedevZhuzhoma2016}, there is a preserving orientation Morse-Smale diffeomorphism $\varphi: N^n\to N^n$ with $NW(\varphi)$ consisting of a unique saddle $\sigma$, a source $\alpha$, and a sink $\omega$. Let us delete small $n$-balls $b_{\alpha}$, $b_{\omega}$ containing $\alpha$ and $\omega$ respectively such that $\sigma\in N^n\setminus(b_{\alpha}\cup b_{\omega})$. Let us identify $\partial b_{\alpha}$ with $S^{n-1}(p^j_t)$ and
$\partial b_{\omega}$ with $S^{n-1}(p^i_s)$. Due to condition (5), the determinants of $A_i$ and $A_j$ are positive. This implies that the restrictions $\varphi_i|_{S^{n-1}(p^i_s)}$, $\varphi_j|_{S^{n-1}(p^j_t)}$ preserve orientation. By construction, the restrictions $\varphi_{\partial b_{\omega}}$, $\varphi_{\partial b_{\alpha}}$ preserve orientation as well. This follows that slightly modifying $\varphi$, $\varphi_i$, and $\varphi_j$, one can get a diffeomorphism
 $$\varphi_{ij}: U(\Omega_i)\cup\left(N^n\setminus(b_{\omega}\cup b_{\alpha}\right)\cup U(\Omega_j)\to U(\Omega_i)\cup\left(N^n\setminus(b_{\omega}\cup b_{\alpha}\right)\cup U(\Omega_j)$$
which is an extension of $\varphi_i$ and $\varphi_j$ such that $\varphi_{ij}(\sigma)=\sigma$. Continuing by this way, one gets a manifold $M^n$ and diffeomorphism  $f\in\mathbb{P}_k(M^n;0,0,1)$ such that $\gamma=\Gamma_{\mathbb{P}}(f)$.

Due to condition (4), the manifold $M^n$ is connected.
The statements concerning a topological structure immediately follows from Lemma \ref{top-struct-for-P}.
$\Box$

\medskip
\textsl{Proof of Theorem \ref{thm:classif-for-M}} is similar to the proof of Theorem \ref{thm:classif-for-P}. Instead of the quality $\widehat{K_{\sigma}}=K_{\sigma}\cup B_1^n\cup B_2^n=N^n$, one uses the quality $\widehat{K_{\sigma}}=K_{\sigma}\cup B_1^n\cup B_2^n\cup B_3^n=\mathbb{S}^n$, and Proposition \ref{prop:sphere-with-holes-homeo} to extend mapping to $\mathbb{S}^n$. We omit details for the Reader.
$\Box$

\bigskip\noindent
\textit{E-mail:} medvedev-1942@mail.ru

\noindent
\textit{E-mail:} zhuzhoma@mail.ru

\end{document}